\documentclass[a4paper, 11 pt]{article}                                                          
\newif\ifarticle
\articletrue

\usepackage[margin=1.3in]{geometry}

\usepackage{epsfig}
\usepackage{booktabs}
\usepackage{amsmath}
\usepackage{amssymb}
\usepackage{enumerate}
\usepackage{url}
\usepackage{caption}
\usepackage{subcaption}

\newtheorem{case}{Case}
\newtheorem{thm}{Theorem}
\newtheorem{fact}{Fact}
\newtheorem{prp}{Proposition}

\newtheorem{ass}{Assumption}
\newtheorem{lem}{Lemma}
\newtheorem{rem}{Remark}
\newtheorem{alg}{Algorithm}
\newtheorem{defin}{Definition}
\newenvironment{pf}{\smallbreak\noindent{\it Proof. }}{\hfill$\Box$\smallbreak}

\makeatletter\@addtoreset{case}{thm}\makeatother
\makeatletter\@addtoreset{case}{subsection}\makeatother

\makeatletter\@addtoreset{subcase}{case}\makeatother

\usepackage{tikz}
\usetikzlibrary{calc}
\usepackage[utf8]{inputenc}
\usepackage{pgfplots}
\usepackage{mathtools}
\usepackage{multicol}
\usepackage{cases}

\newcommand{\fix}{{\rm{fix}}}

\newcommand{\uset}{{\mathcal{U}}}
\newcommand{\vset}{{\mathcal{V}}}
\newcommand{\cset}{{\mathcal{C}}}

\newcommand{\ldef}{\coloneqq}
\newcommand{\rdef}{\eqqcolon}
\newcommand{\re}{\text{\normalfont Re}\,}
\newcommand{\sign}[1]{\text{\normalfont sgn}(#1)}
\newcommand{\cmt}[1]{}
\newcommand{\note}[1]{#1}

\newcommand{\norm}[1]{\left\lVert#1\right\rVert}

\title{\LARGE \bf
Optimal Convergence Rates for Generalized Alternating Projections
}

\author{Mattias F\"{a}lt$^*$ and Pontus Giselsson
\thanks{Department of Automatic Control, Lund University, Sweden. Email:
{\tt\small \{mattiasf,pontusg\}@control.lth.se}
Both authors are financially supported by the Swedish Foundation for Strategic
Research and are members of the LCCC Linneaus Center at Lund University.
}%
}

\begin{document}
\newlength\figureheight
\newlength\figurewidth
\setlength{\figurewidth}{0.86\linewidth}
\setlength{\figureheight }{5cm}
\newcommand{\figurecaptionreduction}{\vspace{0mm}}

\maketitle

\begin{abstract}
Generalized alternating projections is an algorithm that alternates
relaxed projections onto a finite number of sets to find a point in
their intersection. We consider the special case of two linear
subspaces, for which the algorithm reduces to a matrix iteration. For
convergent matrix iterations, the asymptotic rate is linear and decided
by the magnitude of the subdominant eigenvalue. In this paper, we show
how to select the three algorithm parameters to optimize this magnitude,
and hence the asymptotic convergence rate. The obtained rate depends on
the Friedrichs angle between the subspaces and is considerably better
than known rates for other methods such as alternating projections and
Douglas-Rachford splitting. We also present an adaptive scheme that,
online, estimates the Friedrichs angle and updates the algorithm
parameters based on this estimate. A numerical example is provided that
supports our theoretical claims and shows very good performance for the
adaptive method.
\end{abstract}

\section{Introduction}
Many methods for finding a point in the
intersection of a finite number of sets exist.
Notable examples include alternating projections~\cite{vonNeumann,Deutsch1992},
its generalization, generalized alternating projections,
that allows for relaxed projections~\cite{GAP_Agmon, GAP_Motzkin, GAP_Bregman},
Dykstra's algorithm~\cite{Boyle1986}, Douglas-Rachford splitting~\cite{DouglasRachford,LionsMercier1979},
and its dual algorithm ADMM~\cite{Glowinski1975,BoydDistributed}.
Considerable effort has gone into understanding and analyzing
performance and convergence rates of these methods.
Convex and nonconvex feasibility problems have been analyzed in~\cite{Phan20142,hesse2013nonconvex},
and convex optimization and monotone inclusion problems in~\cite{LionsMercier1979,Davis_Yin_2014,gisBoydTAC2014metric_select2,Giselsson2017}.

For feasibility problems with two subspaces,
it has been long known that the standard alternating projection method
converges linearly with exact rate being the squared
Friedrichs angle~\cite{Deutsch1995}.
The Friedrichs angle is the smallest non-zero principal angle between the subspaces,
see~\cite{Deutsch1992} for background on principal angles.
More recently, it was shown in~\cite{Bauschke_lin_rate_Friedrich} that the Douglas-Rachford algorithm
converges with a rate given by the Friedrichs angle.

These projection based algorithms reduce to matrix iterations
when the two sets are subspaces.
This was exploited in~\cite{Bauschke_opt_rate_matr}, where sharp convergence rates for matrices are provided.
They apply their results to find optimal parameters for
the generalized alternating projections method.
Two of the parameters are kept fixed and they optimize over the third.

In this paper, we extend the results of~\cite{Bauschke_opt_rate_matr}.
We optimize the sharp convergence rate for the generalized alternating projection method
over all three algorithm parameters.
The obtained optimal rate turns out to be significantly better
than the ones considered in~\cite{Bauschke_opt_rate_matr}.
The optimal parameters in our setting also depends on the Friedrichs angle.
This angle is of course not known a priori.
Therefore, we have developed an adaptive scheme
that estimates the Friedrichs angle during the course of the iterations.
Under easily achievable assumptions on the starting point of the algorithm,
we show that it is always a conservative estimate of the true Friedrichs angle.
Indeed, in examples we see that the estimated angle approaches the Friedrichs angle.

The intention of this work is not to present a new method for solving
linear systems of equations.
It is rather a starting point to optimize local linear convergence behavior
for the generalized alternating projection method, when solving,
e.g., problems with affine and conic constraints.
Such feasibility problems can solve essentially any convex optimization problem,
by first reformulating the problem as a cone program
(which is done in the CVX modeling languages~\cite{cvx_v3,cvxpy,convexjl}),
and then use primal dual embedding, as in~\cite{SCS}.
The local convergence analysis of such problems is
outside the scope of this paper.
Encouraging results have, however, been presented,
e.g., in~\cite{act_ident_20152} and~\cite{even_lin_conv_20132}.
They show that the local linear convergence rate for
Douglas-Rachford splitting for specific convex optimization problems
is exactly the Friedrichs angle, i.e., the same as for subspaces.
The results rely on sufficient local smoothness or
polyhedral/affine sets and finite identification of active sets or manifolds.
The finite identification property implies that locally,
the problem reduces essentially to an affine subspace intersection problem.

We verify the theoretical results on numerical examples and
demonstrate that the generalized alternating projections with optimal
parameters performs significantly better
than with previously studied parameters in, e.g.,~\cite{Deutsch1992,Bauschke_opt_rate_matr}.
We also observe that the proposed adaptive method performs in line with
the method with optimal parameters.

\section{Preliminaries}

Let the inner product and induced norm be denoted by $\left\langle u,v\right\rangle$  and $\norm{v}\ldef\sqrt{\left<v,v\right>}$ for vectors $u,v\in \mathbb{R}^n$.
Let the set of eigenvalues for a matrix $A\in\mathbb{R}^{n\times n}$ be denoted by $\sigma(A)$,
the spectral radius as $\rho(A)\ldef \max\left\{ \left|\lambda\right|\mid\lambda\in\sigma(A)\right\}$
and let $\norm{A}$ be the operator norm
$\norm{A}\ldef \sup_{x\in\mathbb{R}^n : \norm{x}=1}\norm{Ax}$.
$P_{\cset}$ is the
orthogonal projection onto a closed, convex and nonempty set $\cset$,
i.e. $P_{\cset}x=\text{argmin}_{y\in\cset}\left\{ \norm{x-y}\right\}.$

The following definitions and facts follow closely those in the related work~\cite{Bauschke_opt_rate_matr}.
\begin{defin}\label{def:principal}
The \emph{principal angles} $\theta_{k}\in[0,\pi/2],\,k=1,\dots,p$ between two subspaces
$\uset,\vset\in\mathbb{R}^n$, where $p=\min(\dim\uset,\dim\vset)$, are recursively defined by
\begin{eqnarray*}
\cos\theta_{k}\, & \ldef  & \max_{u_{k}\in\uset,\,v_{k}\in\vset}\left\langle u_{k},v_{k}\right\rangle \\
 & \text{s.t.} & \norm{u_{k}}=\norm{v_{k}}=1,\\
 &  & \left\langle u_{k},v_{i}\right\rangle =\left\langle u_{i},v_{k}\right\rangle =0,\forall\,i=1,...,k-1.
\end{eqnarray*}
\end{defin}
\begin{fact}\label{fct:friedrichs}
\cite[Def 3.1, Prop 3.3]{Bauschke_opt_rate_matr}
The principal angles
are unique and satisfy $0\leq\theta_{1}\leq\theta_{2}\leq\dots\theta_{p}\leq\pi/2$.
The angle $\theta_{F}\ldef \theta_{s+1}$, where $s=\text{dim}(\vset\cap\uset)$, is the \emph{Friedrichs angle} and it
is the smallest non-zero principal angle.
\end{fact}
\begin{defin}
$A\in\mathbb{R}^{n\times n}$ is $\emph{\text{linearly convergent}}$
to $A^{\infty}$ with $\emph{\text{linear convergence rate }}$$\mu\in[0,1)$
if there exist $M,N>0$ such that
\[
\norm{A^{k}-A^{\infty}}\leq M\mu^{k}\quad\forall k>N,\,k\in\mathbb{N}.
\]
\end{defin}
\begin{defin}
\cite[Fact 2.3]{Bauschke_opt_rate_matr}  For $A\in\mathbb{R}^{n\times n}$ we say that
$\lambda\in\sigma(A)$ is \emph{semisimple} if $\text{ker}(A-\lambda I)=\text{ker}(A-\lambda I)^{2}.$
\end{defin}
\begin{fact}\label{fct:limitexists}
\cite[Fact 2.4]{Bauschke_opt_rate_matr} For $A\in\mathbb{R}^{n\times n}$, the limit
$A^{\infty}\ldef \lim_{k\rightarrow\infty}A^{k}$ exists if and only if
\end{fact}
\begin{itemize}
\item $\rho(A)<1$ or
\item $\rho(A)=1$ and $\lambda=1$ is semisimple and the only eigenvalue
on the unit circle.
\end{itemize}
\begin{defin}\label{def:subdominant}
\cite[Def. 2.10]{Bauschke_opt_rate_matr} Let $A\in\mathbb{R}^{n\times n}$ be a (non\-expansive) matrix and define
\[
\gamma(A)\ldef\max\left\{|\lambda|\, \mid\, \lambda\in\{0\}\cup\sigma(A)\setminus\{1\}\right\} .
\]
Then $\lambda\in\sigma(A)$ is a \emph{subdominant eigenvalue} if  $|\lambda|=\gamma(A)$.
\end{defin}%

\begin{fact}
\label{fact:Convergent-at-rate}\cite[Thm. 2.12]{Bauschke_opt_rate_matr}
If $A\in\mathbb{R}^{n\times n}$ is convergent to $A^{\infty}$
then
\end{fact}
\begin{itemize}
\item $A$ is linearly convergent with any rate $\mu\in(\gamma(A),1)$
\item If $A$ is linearly convergent with rate $\mu\in[0,1)$, then $\mu\in[\gamma(A),1)$.
\end{itemize}

%

\section{Optimal parameters for GAP}\label{sec:optimalPars}
Let the relaxed projection onto a set $\mathcal{C}$, with relaxation parameter
$\alpha$, be defined as $P_{\cset}^{\alpha}\ldef (1-\alpha)I+\alpha P_{\cset}$.
The generalized alternating projections (GAP)~\cite{GAPLS} for two
closed, convex and nonempty sets $\uset$
and $\vset$, with $\uset\cap\vset\neq\emptyset$, is then defined by the iteration
\begin{equation}\label{eq:GAP}
x^{k+1}\ldef Sx^{k},
\end{equation}
where
\begin{equation}\label{eq:GAP2}
S=(1-\alpha)I+\alpha P_{\uset}^{\alpha_{2}}P_{\vset}^{\alpha_{1}}=:\,(1-\alpha)I+\alpha T.
\end{equation}

The operator $S$ is averaged and the iterates converge to the fixed-point set $\fix S$
under the following assumption, see e.g.~\cite{GAPLS} where these results are collected.
\begin{ass}\label{ass:alpha}
    Assume that $\alpha\in(0,1]$, $\alpha_{1},\alpha_{2}\in(0,2]$ and that either of the following holds
    \begin{enumerate}[{\it{A}}1.]
    \item $\alpha_{1},\alpha_{2}\in(0,2)$
    \item $\alpha\in(0,1)$ with either $\alpha_1\neq2$ or $\alpha_2\neq2$
    \item $\alpha\in(0,1)$ and $\alpha_{1}=\alpha_{2}=2$
    \end{enumerate}
\end{ass}

To study the convergence rate of $S$,
and its dependence on the parameters $\alpha_1, \alpha_2$ and $\alpha$,
we need to characterize the eigenvalues of $S$.
To this end, we state the following proposition, as found in~\cite[Prop. 3.4]{Bauschke_opt_rate_matr}.
    \begin{prp}\label{prp:projections}
        Let $\uset$ and $\vset$ be affine subspaces in $\mathbb{R}^{n}$
        satisfying $p\ldef \dim(\uset)$, $q\ldef \dim(\vset)$,
        where $p\leq q$, $p+q< n$ and $p,q\geq1$.
        Then, the projection matrices $P_{\uset}$
        and $P_{\vset}$ become

    \begin{align}
    P_{\uset} & =D\begin{pmatrix}I_{p} & 0 & 0 & 0\\
    0 & 0_{p} & 0 & 0\\
    0 & 0 & 0_{q-p} & 0\\
    0 & 0 & 0 & 0_{n-p-q}
    \end{pmatrix}D^{*},\label{eq:Pu}\\
    P_{\vset} & =D\begin{pmatrix}C^{2} & CS & 0 & 0\\
    CS & S^{2} & 0 & 0\\
    0 & 0 & I_{q-p} & 0\\
    0 & 0 & 0 & 0_{n-p-q}
    \end{pmatrix}D^{*}\label{eq:Pv}
    \end{align}
    and
    \begin{equation}
    P_{\uset}P_{\vset}=D\begin{pmatrix}C^{2} & CS & 0 & 0\\
    0 & 0_{p} & 0 & 0\\
    0 & 0 & 0_{q-p} & 0\\
    0 & 0 & 0 & 0_{n-p-q}
    \end{pmatrix}D^{*},
    \end{equation}
    where $C$ and $S$ are diagonal matrices containing the cosine and sine
    of the principal angles $\theta_{i}$, i.e.
    \begin{align*}
    S&=\text{\normalfont diag}(\sin\theta_{1},\dots,\sin\theta_{p}),\\
    C&=\text{\normalfont diag}(\cos\theta_{1},\dots,\cos\theta_{p}),
    \end{align*}
    and $D\in\mathbb{R}^{n\times n}$ is an orthogonal matrix.%
\end{prp}

Under the assumptions in Proposition~\ref{prp:projections},
the linear operator $T$, implicitly defined in~\eqref{eq:GAP2}, becomes
\begin{eqnarray*}
T & = & P_{\uset}^{\alpha_{2}}P_{\vset}^{\alpha_{1}}=((1-\alpha_{2})I+\alpha_{2}P_{\uset})((1-\alpha_{1})I+\alpha_{1}P_{\vset})\\
 & = & (1-\alpha_{2})(1-\alpha_{1})I+\alpha_{2}(1-\alpha_{1})P_{\uset}\\
 &  & +\alpha_{1}(1-\alpha_{2})P_{\vset}+\alpha_{1}\alpha_{2}P_{\uset}P_{\vset}\\
 & = & D\,\text{blkdiag}(T_{1},T_{2},T_{3})\,D^{*}
\end{eqnarray*}
where
\begin{align}
T_{1} & =\begin{pmatrix}I_{p}-\alpha_{1}S^{2} & \alpha_{1}CS\\
\alpha_{1}(1-\alpha_{2})CS & (1-\alpha_{2})(I_{p}-\alpha_{1}C^{2})
\end{pmatrix},\label{eq:allblocks}\\
T_{2} & =(1-\alpha_{2})I_{q-p},\quad T_{3}=(1-\alpha_{2})(1-\alpha_{1})I_{n-p-q}.\nonumber
\end{align}
The rows and columns of $T_1$ can be reordered so that it is a block-diagonal
matrix with blocks
\begin{equation}
T_1^i=\begin{pmatrix}1-\alpha_{1}s_{i}^{2} & \alpha_{1}c_{i}s_{i}\\
\alpha_{1}(1-\alpha_{2})c_{i}s_{i} & (1-\alpha_{2})(1-\alpha_{1}c_{i}^{2})
\end{pmatrix},\quad i\in1,\dots,p\label{eq:T1matrix}
\end{equation}
where $s_{i}\ldef \sin\theta_{i},\,c_{i}\ldef \cos\theta_{i}$.
The eigenvalues of $T$ are therefore
    $\lambda^{3}\ldef (1-\alpha_{2})$,
    $\lambda^{4}\ldef (1-\alpha_{2})(1-\alpha_{1})$,
and for every $T_1^1$
\begin{alignat}{1}
\lambda_{i}^{1,2} & =\frac{1}{2}\left(2-\alpha_{1}-\alpha_{2}+\alpha_{1}\alpha_{2}c_{i}^{2}\right)\label{eq:eig12}\\
 & \,\,\,\pm\sqrt{\frac{1}{4}\left(2-\alpha_{1}-\alpha_{2}+\alpha_{1}\alpha_{2}c_{i}^{2}\right)^{2}-(1-\alpha_{1})(1-\alpha_{2})}.\nonumber
\end{alignat}

\begin{rem}\label{rem:a1a2}
    The property $p\leq q$ was used to arrive at these results. If instead $p > q$,
    we reverse the definitions of $P_\uset$ and $P_\vset$ in Proposition~\ref{prp:projections}.
    Noting that $\sigma(T)=\sigma(T^T)$,
    we get a new block-diagonal matrix $\bar T$ with blocks
    $\bar{T}_1=T_1^T$, $\bar{T}_3=T_3^T$ and $\bar{T}_{2} = (1-\alpha_{1})I_{p-q}$.
    Therefore, the matrix will have eigenvalues in either
    $1-\alpha_1$ or $1-\alpha_2$ depending on the dimensions of $\uset$ and $\vset$.
\end{rem}

Motivated by Fact~\ref{fact:Convergent-at-rate},
we are looking for parameters that minimize the magnitude of the subdominant eigenvalues.
We will do this for both cases in Remark~\ref{rem:a1a2}.
In the following sequence of theorems,
we will show that the optimal parameters are
%
\begin{align}\label{eq:optpar}
    \alpha=1,\quad\alpha_{1}=\alpha_{2}=\alpha^{*}\ldef \frac{2}{1+\sin{\theta_{F}}},
\end{align}
and that the subdominant eigenvalues have magnitude $\gamma(S)=\gamma^{*}$, where
\begin{align}\label{eq:optrate}
\gamma^* \ldef \frac{1-\sin\theta_{F}}{1+\sin\theta_{F}}.
\end{align}
\begin{thm}\label{thm:gapP1-1}
The GAP operator $S$ in~\eqref{eq:GAP2}
with $\alpha,\alpha_1,\alpha_2$ as defined in~\eqref{eq:optpar}
satisfies $\gamma(S)=\gamma^{*}$
and is linearly convergent with any rate $\mu\in\left(\gamma^{*},1\right)$.
\end{thm}
A proof is located in Appendix~\ref{app:gapP1-1}.

We now show that no other choices of $\alpha,\alpha_{1},\alpha_{2}$
can achieve a lower linear convergence rate under the assumption that the
relative dimension of $\uset$ and $\vset$ is unknown.
Motivated by this, we formulate the following assumption.
\begin{ass}\label{ass:dimUV}
    Suppose that $\uset$ and $\vset$ are linear subspaces and that the dimensions
    $p\ldef\dim(\uset),\, q\ldef\dim(\vset)$ satisfy
    $p,q\in\{1,\dots,n-1\}$ and consider the cases:
    \begin{align*}
    \text{B1: } p< q, \quad \text{B2: } p=q, \,\,\,  \text{and} \,\,\,\, \text{B3: } p>  q.
    \end{align*}
\end{ass}

\begin{prp}\label{prp:eigset}
    To optimize the convergence rate of $S$, for all cases in Assumption~\ref{ass:dimUV},
    it is necessary to minimize the largest modulus of the eigenvalues in the set
    \begin{align}\label{eq:fullset}
    \left(\{ \lambda_{i}^{1,2}\} _{i\in1,\dots,p}\cap\left\{ 1-\alpha_{2}, 1-\alpha_{1},\,(1-\alpha_{2})(1-\alpha_{1})\right\}\right)\setminus\{1\}.
    \end{align}
\end{prp}
\begin{pf}
    These are the eigenvalues from the matrices in~\eqref{eq:allblocks} together with $1-\alpha_1$,
    as motivated in Remark~\ref{rem:a1a2}.
    If we let $\gamma_1=\gamma(S)$ under assumption B1, $\gamma_2=\gamma(S)$ under B2,
    and $\gamma_3=\gamma(S)$ under B3,
    it follows, from Remark~\ref{rem:a1a2}, that the largest modulus of the eigenvalues in~\eqref{eq:fullset} is
    equal to $\max(\gamma_1,\gamma_2,\gamma_3)$.
\end{pf}
Next, we show that the rate obtained in Theorem~\ref{thm:gapP1-1} is indeed optimal.
\begin{thm}\label{thm:gapiff}
The GAP operator $S$ in~\eqref{eq:GAP2} with $\theta_F<\pi/2$ and $\alpha_1,\alpha_2,\alpha>0$
is linearly convergent with any rate $\mu\in\left(\gamma^{*},1\right)$, for all cases in Assumption 2,
if and only if $\alpha,\alpha_{1},\alpha_{2}$ are chosen as in~\eqref{eq:optpar}.

\end{thm}
A proof is located in Appendix~\ref{app:gapiff}.
\begin{rem}
    The case with $\theta_F=\pi/2$ is trivial and results in convergence in one iteration with the optimal parameters.
    This case is excluded from the theorem since there are also other methods that achieve the same rate.
    We also exclude the cases when either of $\alpha_1,\alpha_2,\alpha$ are non-positive, since
    such choices typically result in a non-convergent algorithm.
    The assumption on the parameters is, however, less restrictive than Assumption~\ref{ass:alpha}.
\end{rem}
\begin{rem}
    The result is derived under the assumption that both $1-\alpha_{2}$ and $1-\alpha_{1}$
    are considered, i.e. $q < p$ and $q > p$ respectively (see Remark~\ref{rem:a1a2}).
    The same result follows in either of these cases if we instead assume that $\theta_p=\pi/2$,
    which is a safe assumption if we do not know the largest principal angle.
\end{rem}
We now state the convergence rate of the sequence $x^k$.
\begin{thm}\label{thm:linconv}
The sequence $x^{k+1}\ldef S x^k$ with optimal parameters
$\alpha=1, \alpha_1=\alpha_2=\frac{2}{1+\sin\theta_F}$
converges linearly to $x^{*}\ldef P_{\fix S}x^{0}$ according to
\begin{align}
    \norm{x^{k}-x^{*}} \leq \mu^k\norm{x^0}\quad \forall k\geq N,
\end{align}
with any rate $\mu\in(\gamma^*,1)$, for $\gamma^*$ in~\eqref{eq:optrate},
i.e., $x^k$ is R-linearly convergent to $x^*$.
\end{thm}
A proof is located in Appendix~\ref{app:linconv}.

\begin{rem}
    For linear subspaces~$\uset,\vset$,
    under the Assumption~\ref{ass:alpha} case A1 or A2,
    we have $\fix S=\uset\cap\vset$, see e.g.~\cite{GAPLS}.
      For case A3 we have $\fix S=\vset\cap\uset + (\vset^\perp\cap\uset^\perp)$,
    see~\cite{Bauschke_lin_rate_Friedrich}.
\end{rem}

\section{Comparison with other choices of parameters}\label{sec:otherchoises}
In Section~\ref{sec:optimalPars}, we derive, for two linear subspaces,
the optimal parameters for the
generalized alternating projections method.
These parameters are optimal under the assumption that the relative dimensions of the
two subspaces are unknown, or that the largest principal angle $\theta_p=\pi/2$.
There are other methods that can perform better if these assumptions are not true.
For example, if $\dim\uset\leq\dim\vset$, the parameters
\begin{equation}\label{eq:GAP2a}
    \alpha=1,\quad \alpha_1=2,\quad \alpha_2=\frac{2}{1+\sin(2\theta_F)},
\end{equation}
(referred to as GAP2$\alpha$ in Section~\ref{sec:numerical}) result in that most eigenvalues have modulus
\begin{equation}
    \frac{\cos\theta_F-\sin\theta_F}{\cos\theta_F+\sin\theta_F}.
\end{equation}
This rate is better than $\gamma^*$, although marginally for small $\theta_F$.
However, if the largest principal angle, $\theta_p$, is large enough,
the corresponding eigenvalues will approach $-1$. This choice will
then converge much slower than the optimal method in Section~\ref{sec:optimalPars}.
This is observed in the numerical example in Section~\ref{sec:numerical}.

When $\dim\uset\leq\dim\vset$,
it is sometimes possible to get even better performance by
selecting $\alpha_2>2$.
However, this method is not convergent if $\dim\uset>\dim\vset$,
and it would generally not be convergent for general convex sets.

\newcommand{\DrawEigs}[5]{
  \pgfmathsetmacro{\phif}{#4}
  \draw[->] (0,-1.1) -- (0,1.1);
  \node at (0.2,1.15) {Im};
  \draw[->] (-1.1,0) -- (1.1,0);
  \node at (1.2,0.15) {Re};
  \draw[thick] (0,0) circle (1cm);
  \draw (0,0) circle (1cm);
  \pgfmathsetmacro{\N}{30}
  \pgfmathsetmacro{\a}{#1}
  \pgfmathsetmacro{\aa}{#2}
  \pgfmathsetmacro{\ab}{#3}
  \foreach \x in {0,...,\N}{
    \pgfmathsetmacro{\t}{(\phif*(\N-\x) + #5*\x)/(\N)}
    \pgfmathsetmacro{\px}{1/2*(2-\aa-\ab+\aa*\ab*cos(\t)*cos(\t))}
    \pgfmathsetmacro{\py}{\px*\px - (1-\aa)*(1-\ab)}
    \pgfmathsetmacro{\rea}{\ifdim \py pt>0pt \px + sqrt(\py) \else \px \fi }
    \pgfmathsetmacro{\ima}{\ifdim \py pt>0pt 0.0 \else sqrt(-\py) \fi }
    \pgfmathsetmacro{\reb}{\ifdim \py pt>0pt \px - sqrt(\py) \else \px \fi }
    \pgfmathsetmacro{\imb}{\ifdim \py pt>0pt 0.0 \else -sqrt(-\py) \fi }
    \node[circle, fill, inner sep=1pt, blue!\t!red] at ({\a*(\rea-1)+1}, {\a*\ima}) {};
    \node[circle, fill, inner sep=1pt, blue!\t!red] at ({\a*(\reb-1)+1}, {\a*\imb}) {};
  }
    \node[green, circle, fill, inner sep=1pt] at ({\a*(1-\ab)+(1-\a)}, 0) {};
    \node[green, circle, fill, inner sep=1pt] at ({\a*(1-\ab)*(1-\aa)+(1-\a)}, 0) {};
}

\pgfmathsetmacro{\tf}{8.195}
\pgfmathsetmacro{\tbig}{90}
\ifarticle
    \newcommand{\figscale}{1.95}
    \newcommand{\subfigwidth}{0.32\textwidth}
\else
    \newcommand{\figscale}{1.9}
    \newcommand{\subfigwidth}{0.24\textwidth}
\fi
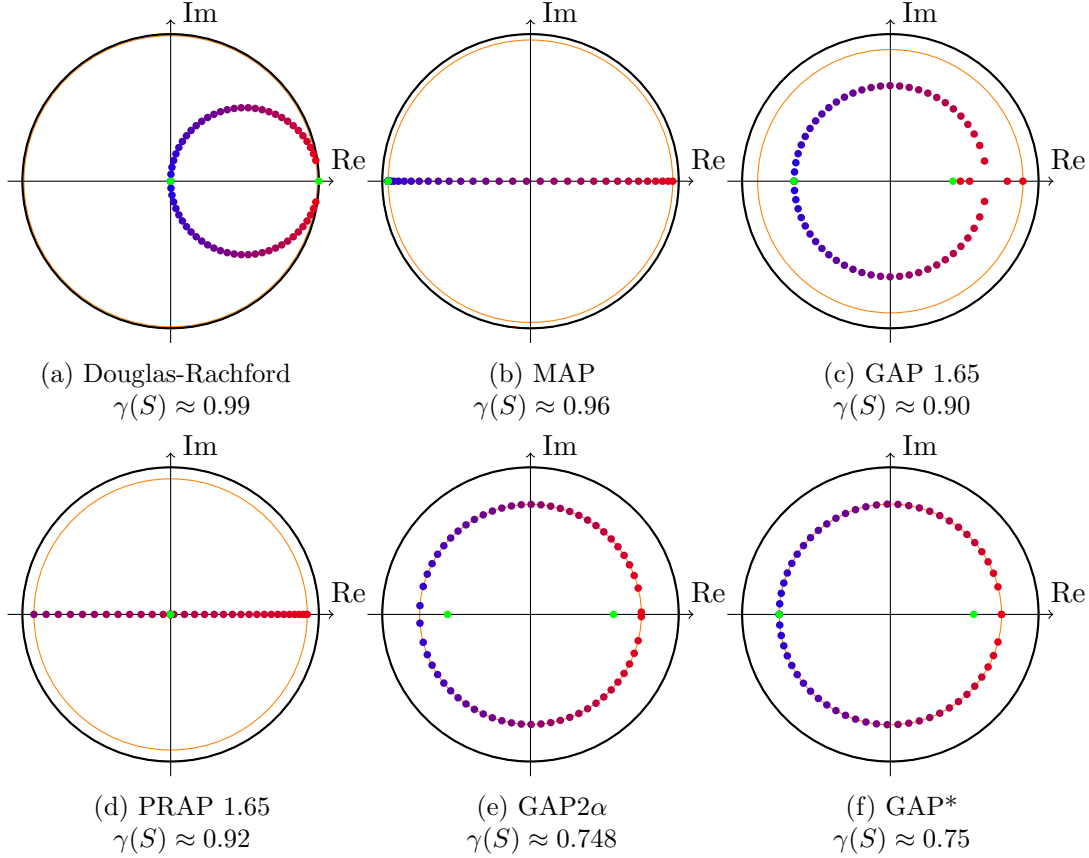
\begin{figure*}
    \tikzset{every picture/.style={scale=\figscale, cross/.style={path picture={
      \draw
      (path picture bounding box.south east) -- (path picture bounding box.north west) (path picture bounding box.south west) -- (path picture bounding box.north east);
    }}}}
    \centering
    \begin{subfigure}{\subfigwidth}
        \begin{tikzpicture}
            \pgfmathsetmacro{\r}{cos(\tf)}
            \draw[orange] (0,0) circle (\r cm);
            \DrawEigs{0.5}{2.}{2.}{\tf}{\tbig}
        \end{tikzpicture}
        \caption{Douglas-Rachford\\\centering{$\gamma(S)\approx0.99$}}
    \end{subfigure}
    \begin{subfigure}{\subfigwidth}
        \begin{tikzpicture}
            \pgfmathsetmacro{\r}{(1-sin(\tf)*sin(\tf))/(1+sin(\tf)*sin(\tf))}
            \draw[orange] (0,0) circle (\r cm);
            \DrawEigs{{2/(1+sin(\tf)*sin(\tf))}}{1.}{1.}{\tf}{\tbig}
        \end{tikzpicture}
        \caption{MAP\\$\gamma(S)\approx0.96$}
    \end{subfigure}
    \begin{subfigure}{\subfigwidth}
        \begin{tikzpicture}
            \pgfmathsetmacro{\tbigg}{45}
            \pgfmathsetmacro{\r}{0.8952}
            \draw[orange] (0,0) circle (\r cm);
            \DrawEigs{1}{1.65}{1.65}{\tf}{\tbig}
        \end{tikzpicture}
        \caption{GAP 1.65\\$\gamma(S)\approx0.90$}
    \end{subfigure}
    \ifarticle
        \begin{subfigure}{\subfigwidth}
            \begin{tikzpicture}
                \pgfmathsetmacro{\tbigg}{45}
                \pgfmathsetmacro{\r}{(sin(\tbigg)*sin(\tbigg)-sin(\tf)*sin(\tf))/(sin(\tbigg)*sin(\tbigg)+sin(\tf)*sin(\tf))}
                \draw[orange] (0,0) circle (\r cm);
                \DrawEigs{1}{{2/(sin(\tbigg)*sin(\tbigg)+sin(\tf)*sin(\tf))}}{1.}{\tf}{\tbigg}
            \end{tikzpicture}
            \caption{PRAP 1.65\\$\gamma(S)\approx0.92$}
        \end{subfigure}
        \begin{subfigure}{\subfigwidth}
            \begin{tikzpicture}
                \pgfmathsetmacro{\tbigg}{81.5}
                \pgfmathsetmacro{\r}{(cos(\tf)-sin(\tf))/(cos(\tf)+sin(\tf))}
                \draw[orange] (0,0) circle (\r cm);
                \DrawEigs{1}{2}{{2/(1+sin(2*\tf)}}{\tf}{\tbigg}
            \end{tikzpicture}
            \caption{GAP2$\alpha$\\$\gamma(S)\approx0.748$}
        \end{subfigure}
    \fi
    \begin{subfigure}{\subfigwidth}
        \begin{tikzpicture}
            \pgfmathsetmacro{\r}{(1-sin(\tf))/(1+sin(\tf))}
            \draw[orange] (0,0) circle (\r cm);
            \DrawEigs{1.}{1.75}{1.75}{\tf}{\tbig}
        \end{tikzpicture}
        \caption{GAP*\\$\gamma(S)\approx0.75$}
    \end{subfigure}
    %
    %
    %
    %
    \caption{\label{fig:eigs}
    Convergence rates for different methods, as described in Section~\ref{sec:otherchoises},
    for $\theta_F\approx0.14\,(8.8^\circ)$.
    The eigenvalues corresponding to the principal angels are shown for $30$ angles,
    evenly spaced from $\theta_F$ to $\pi/2$, as dots from red to blue.
    The eigenvalues corresponding to $(1-\alpha_2)$ and $(1-\alpha_2)(1-\alpha_1)$ are shown as green dots.
    The radius $\gamma(S)$ is shown in orange.
    GAP1.65 represents GAP with $\alpha=1$ and $\alpha_1=\alpha_2=1.65<\alpha^*=1.75$.
    The partial relaxed alternating projections (PRAP) from Equation~\eqref{eq:prap}, the best algorithm in the previous work~\cite{Bauschke_opt_rate_matr}, is shown under the assumption $\theta_p=\pi/4$.
    We see that the optimal parameters gives a much better result than the previously suggested methods.
    This is achieved by placing the eigenvalues at the same radius.
    Increasing the parameters from the optimal ($\alpha_1=\alpha_2>\alpha^*=1.75$),
    increases the radius of the eigenvalues corresponding to the principal angles.
    If decreased, the result looks like GAP 1.65, where one of the eigenvalues corresponding to
    $\theta_F$ is subdominant.
    GAP2$\alpha$ (Equation~\eqref{eq:GAP2a}) is shown under the assumption $\theta_p\approx0.91\pi/2$. Although it performs slightly better than GAP* under this assumtion, it gets considerably worse if $\theta_p$ increases.
    }
\end{figure*}

In~\cite{Bauschke_opt_rate_matr},
optimal parameters are found by keeping two of the parameters
fixed and optimizing over the third.

The first method is the relaxed alternating projections ($\alpha_1=\alpha_2=1$), which is shown to be optimal for
$\alpha=\frac{2}{1+\sin^2\theta_F}$ with rate
$\gamma=(1-\sin^2\theta_F)/(1+\sin^2\theta_F)$.
This is better than the alternating projections with $\alpha=1$
which is convergent with rate $\gamma=\cos^2\theta_F$~\cite{Deutsch1995}.

The generalized Douglas-Rachford ($\alpha_1=\alpha_2=2$), is shown to be optimal for
$\alpha=0.5$ with rate
$\gamma=\cos\theta_F$.

These rates are considerably worse than the optimal rates, as seen in Figure~\ref{fig:eigs},
especially for small $\theta_F$. The methods are referred to as
MAP and DR in the numerical example in Section~\ref{sec:numerical}.

The partial relaxed alternating projections ($\alpha=\alpha_2=1$) was
was shown to be optimal for
\begin{equation}\label{eq:prap}
    \alpha_2=\frac{2}{\sin^2\theta_p+\sin^2\theta_F},\quad\text{with rate}\quad
    \gamma=\frac{\sin^2\theta_p-\sin^2\theta_F}{\sin^2\theta_p+\sin^2\theta_F}.
\end{equation}
%
This rate is sometimes better than $\gamma^*$ if $\theta_p<\pi/2$,
but not for small enough $\theta_F$. In fact, it is only better if $\sin^2\theta_p<\sin\theta_F$.
It also requires knowledge of $\theta_p$, and is not generally convergent if $\dim\uset>\dim\vset$.

An illustration of where the eigenvalues are located for
\ifarticle
these
\else
a few different
\fi
methods is shown in Figure~\ref{fig:eigs}.

\section{Adaptive generalized alternating projections\label{sec:GAPA}}

The generalized alternating projections algorithm with $\alpha_{1}=1$,
$\alpha_{1}=\alpha_{2}=\frac{2}{1+\sin\theta_{F}}$ is
optimal under the assumption
that the relative dimensions between $\uset$ and $\vset$ is unknown.
However, this parameter choice requires that the Friedrichs angle is known.
This is typically not the case.
In this section, we present an adaptive method that continuously
tries to estimate the Friedrichs angle $\theta_{F}$ and updates $\alpha_{1}$ and $\alpha_{2}$, based on this estimate.

Consider the following estimate of the Friedrichs angle at iteration $k$
\begin{align}\label{eq:est}
    \cos\hat{\theta}^k & \ldef \frac{|\left< x^k-y^k, z^k-y^k\right>|}{\norm{x^k-y^k}\norm{z^k-y^k}},
\end{align}
where $y^k=P^{\alpha_1}_\vset x^k$ and $z^k=P_\uset P_\vset^{\alpha_1}x^k$.
If $x^k=y^k$ or $z^k=y^k$ we define the estimate as $\cos\theta_k\ldef0$.
The estimate is illustrated in Figure~\ref{fig:est}\\
Next, we show that this value is always an overestimation of the Friedrichs angle,
provided that the first iterate is in $\uset+\vset$.
\begin{figure}
    \centering
    \begin{tikzpicture}[dot/.style={circle,inner sep=1pt,fill},
  extended line/.style={shorten >=-#1,shorten <=-#1},
  extended line/.default=1cm,
  every node/.style={inner sep=0pt, outer sep=0pt}, scale=0.8]
\usetikzlibrary{calc}

\node (o) at (0,0) {};
\node[label=right:{$\mathcal{V}$}]  (Vend) at (4,0) {};
\node[label=right:{$\mathcal{U}$}]  (Uend) at (4,2.3) {};

\draw(o) -- (Vend);
\draw (o) -- (Uend);

\node[dot, label=above right:{$x^k$}] (xk) at (3.5,1.4) {$$};
\node[dot, label=above right:{$y_k\ldef P^{\alpha_1}_\mathcal{V}x^k$}] (yk) at (3.5,-1.3) {};
\node[dot, label=above left:{$z_k\ldef P_\mathcal{U}P^{\alpha_1}_\mathcal{V}x^k$}] (zb) at  ($(o)!(yk)!(Uend)$) {};

\draw (xk) -- (yk) {};

\node[dot, label=below right:{}] (yb) at ($(o)!(xk)!(Vend)$) {};

\draw [] ($(o)!(yk)!(Uend)$) -- (yk);

\begin{scope}
\path[clip] (zb) -- (3.5,-1.3) -- (xk);
\fill[red, opacity=0.5, draw=black] (yk) circle (7mm);
\node at ($(yk)+(102:10mm)$) {$\hat\theta^k$};
\end{scope}



\end{tikzpicture}
    \caption{Illustration of the estimate $\hat{\theta}^k$.\label{fig:est}}
\end{figure}
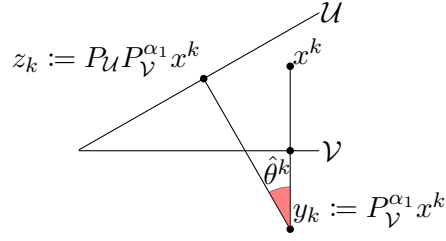
\begin{thm}
    The estimate $\hat{\theta}^k$ in Equation~\eqref{eq:est} always satisfies
    $\hat{\theta}^k \geq \theta_F$ if the starting point $x^0\in \uset + \vset$.
\end{thm}
\begin{pf}
    Assume that $x^k\in \uset + \vset$.
    Since for a projection  it holds that $P_\vset x^k\in \vset$,
    it follows that $y^k = P_\vset^{\alpha_1} x^k$,
    a linear combination of $x^k$ and $P_\vset x^k$, satisfies $y^k\in\uset + \vset$.
    In the same way it follows that $z^k\in \uset + \vset$ and $x^{k+1}\in\uset + \vset$.
    By induction, this must hold for all iterations since $x^0\in \uset + \vset$.

    Let $v_1\ldef x^k-y^k$ and $v_2\ldef z^k-y^k$.
    We have $v_1=x^k-P^{\alpha_1}_\vset x^k=\alpha_1(I-P_\vset)x^k=\alpha_1 P_{\vset^\perp}x^k\in \vset^\perp$
    and in the same way $v_2\in\uset^\perp$.
    We also see that $v_1,v_2\in\uset+\vset$, since they are linear combinations of elements in
    $\uset+\vset$. Noting that
    $\uset+\vset=(\uset^\perp\cap\vset^\perp)^\perp$~\cite[Lem. 2.11]{Deutsch1995} we get,
    \begin{align*}
        v_1\in\uset^\perp\cap(\uset^\perp\cap\vset^\perp)^\perp
        ,\, v_2\in\vset^\perp\cap(\uset^\perp\cap\vset^\perp)^\perp.
    \end{align*}
    Using the definition of the cosine of the Friedrichs angle between two sets $\uset,\vset$~\cite[Def. 2.1]{Deutsch1995}:
    \begin{align*}
      c_F(\uset,\vset) \ldef \max \left\{\frac{|\left<v,u\right>|}{\norm{v}\norm{u}}
      \,:\quad
      \begin{matrix}
        v\in\uset\cap(\uset\cap\vset)^\perp\\
        u\in\vset\cap(\uset\cap\vset)^\perp
      \end{matrix}
      \right\}
    \end{align*}
    and the property $c_F(\uset,\vset)=c_F(\uset^\perp,\vset^\perp)$~\cite[Thm. 2.16]{Deutsch1995} we immediately get
    \begin{align*}
      \cos\hat{\theta}^k = \frac{|\left<v_1,v_2\right>|}{\norm{v_1}\norm{v_2}}
      \leq c_F(\uset^\perp, \vset^\perp) = c_F(\uset, \vset) = \cos\theta_F
    \end{align*}
    where we let $\frac{|\left<v_1,v_2\right>|}{\norm{v_1}\norm{v_2}}\ldef0$
    if $\norm{v_1}=0$ or $\norm{v_2}=0$.\\
    We therefore conclude that
    $\hat{\theta}^k \geq \theta_F$.
\end{pf}

Next, we propose an adaptive version of the generalized alternating projections method:
\begin{alg}\label{alg:adaptivegap}
Let $k=0$, $x^0\in\mathbb{R}^n$ and $\alpha^0\in(0,2)$.
\begin{align*}
    y^{k}                              & \ldef P^{\alpha^k}_\vset x^k\\
    x^{k+1}                              & \ldef P_\mathcal{U}^{\alpha^k}y^{k}\\
    \hat{\theta}^k                     & \ldef \text{\normalfont acos}\frac{|\left< x^k-y^k, x^{k+1}-y^k\right>|}{\norm{x^k-y^k}\norm{x^{k+1}-y^k}}\\
    \alpha^{k+1}                       & \ldef \frac{2}{1+\sin\hat{\theta}^k}
\end{align*}
\end{alg}

We now motivate, without proof,
that the estimate will tend toward $\theta_F$ if $x^0\in\uset+\vset$.

Let $\hat{\theta}^k$ be the current estimate of $\theta_{F}$
and $\alpha_{1}=\alpha_{2}=\frac{2}{1+\sin\hat{\theta}^k}$.
Since $\hat{\theta}^k\geq\theta_{F}$, we get $\alpha_{1}=\alpha_{2}\leq\alpha^*$.
As seen in Figure~\ref{fig:eigs}(c), eigenvalues corresponding to
large principal angles
have radius smaller than $\alpha^{*}-1$.
However smaller principal angles
will have one positive real eigenvalue, and the largest eigenvalue corresponds
to $\theta_{F}$ with real part greater than $\alpha^*-1$. Iterating
the operator should therefore result in convergence to the subspace
spanned by the eigenvectors corresponding to $\theta_{F}$, and the
estimated angle will decrease towards $\theta_{F}$.
This behavior was observed in the numerical example in Section~\ref{sec:numerical}.

We now show that Algorithm~\ref{alg:adaptivegap} is always convergent,
for general convex sets, if it is modified so that $\alpha^k\neq 2$.
This is true if $\hat{\theta}_F>0$ or if the algorithm is modified, for example as
\begin{align*}
  \alpha^k\gets \min\left\{\frac{2}{1+\sin\hat{\theta}^k}, 2-\epsilon\right\},
\end{align*}
for some $\epsilon>0$.
\begin{thm}
  Consider Algorithm~\ref{alg:adaptivegap} for two non-empty, closed, convex sets
  $\uset,\vset$ with $\uset\cap\vset\neq\emptyset$.
  If $\hat{\theta}^k$ satisfies
    $\hat{\theta}^k>0$ for all $k\geq0$ then
  $x^k\rightarrow x^*$ for some $x^*\in\uset\cap\vset$.
\end{thm}
\begin{pf}
  If $\hat{\theta}^k>0$, then $\alpha^{k+1}\neq 2$.
  Thus $\alpha^{k+1}\in(0,2)$ and each iteration is the result of an averaged mapping $S^k$
  with fixed points $\uset\cap\vset$.
  It follows that the iterates converge to the fixed point set $\uset\cap\vset$, see e.g.~\cite{GAPLS}.
\end{pf}

\section{Numerical Example}\label{sec:numerical}

\begin{figure*}
\begin{center}
\def\factor{.23}%
\begin{tikzpicture}
\begin{loglogaxis}[clip mode=individual,
legend entries={GAP*, GAPA, GAP2$\alpha$, DR, MAP, GAP$1.8$},
xmode=log,ymode=log,xmin=0.0004,ymin=2,xmax=1.2,ymax=250000,
width=\textwidth,
xlabel={Friedrichs angle $\theta_F$}, ylabel={Iterations},
y label style={at={(axis description cs:0.02,.5)},anchor=south},
legend cell align=left]

\definecolor{bblue}{rgb}{0.3,0.3,1.0}
\definecolor{lblue}{rgb}{0.45,0.45,1.0}
\definecolor{ggreen}{rgb}{0.0,1.0,0.0}
\definecolor{dgreen}{rgb}{0.0,0.75,0.0}
\definecolor{rred}{rgb}{1.0,0.0,0.0}
\definecolor{lred}{rgb}{1.0,0.2,0.2}
\definecolor{cya}{rgb}{0.0,0.8,0.8}
\definecolor{oorange}{rgb}{1.0,0.7,0.0}
\definecolor{ppurple}{rgb}{0.3,0.0,0.6}
\definecolor{lpurple}{rgb}{0.6,0.0,1.0}

\addlegendimage{only marks, bblue}
\addlegendimage{only marks, oorange}
\addlegendimage{only marks, rred}
\addlegendimage{only marks, dgreen}
\addlegendimage{only marks, ppurple}
\addlegendimage{only marks, cya}
\addplot graphics[xmin=0.0004,ymin=2,xmax=1.2,ymax=250000]{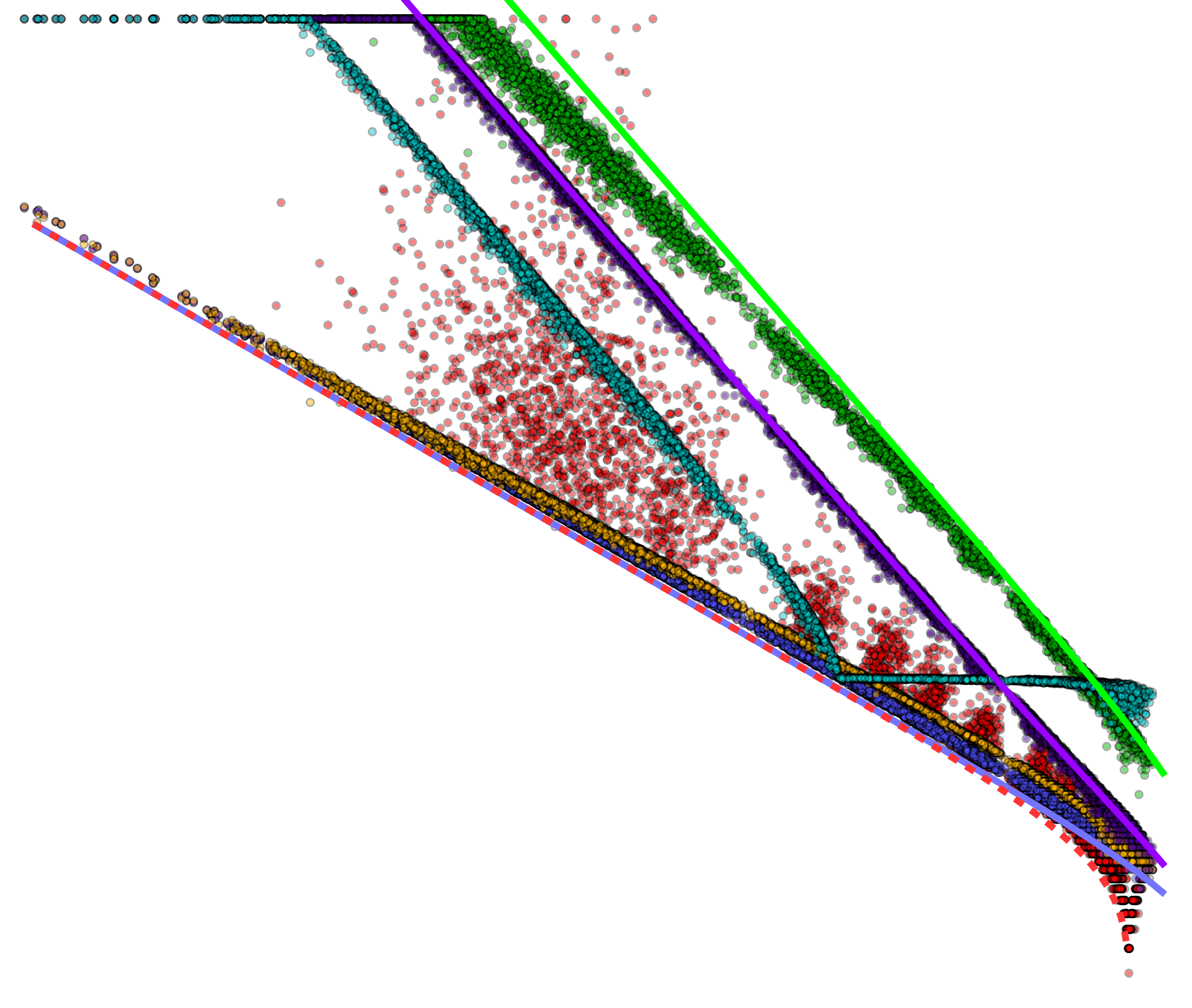};

\end{loglogaxis}

\end{tikzpicture}%
\end{center}
\caption{\label{fig:iters-1-1}Number of iterations for different methods,
as described in Section~\ref{sec:numerical},
plotted against the Friedrichs angle $\theta_{F}$.
The theoretical rates are plotted in lines as the solution to $\gamma(S)^n=10^{-8}$
for GAP*, DR, and MAP. For GAP2$\alpha$ we show the rate (in dashed red line) assuming that
$\theta_p$ is sufficiently small,
according to the discussion in Section~\ref{sec:otherchoises}.
We see that this method can perform better than GAP*,
particularly for large $\theta_F$.
However, since $\theta_p$
is unknown, convergence is sometimes extremely slow.
The convergence for GAP1.8 is constant for small $\theta_F$, but
the convergence rate slows down considerably when
$\theta_F$ decreases to the point where $1.8<\alpha^*$.
We see that GAP* performs in line with the theoretical result,
and considerably better than both DR and MAP.
The adaptive method (GAPA) performs marginally worse than GAP*
for large $\theta_F$.
No difference in the number of iterations can be seen between GAP* and GAPA when $\theta_F$ is small.
}

\end{figure*}

In this section, we compare the theoretical results to numerical experiments.
We have generated a set of problems of the form
\[
\vset=\left\{ x\mid Ax=0\right\} ,\,\uset=\left\{ x\mid Bx=0\right\}
\]
with $A\in\mathbb{R}^{n\times200},\,B\in\mathbb{R}^{100\times200}$.
The matrices are generated with independent normal distributed
elements, with zero mean and unit variance. The initial point $x^{0}$
is randomly chosen in the same way. The dimension of $A$ is selected from
$13$ different categories with $n\in\left\{ 1,\dots,99\right\} $,
and at least $500$ problems are generated for each category, resulting
in over $8000$ different problems. The problems have
Friedrichs angles in the range $\theta_{F}\in(5\cdot10^{-4},1).$

We solve the problem of finding $x\in\uset\cap\vset$ using the following
algorithms:
\begin{itemize}
\item Method of alternating projections (MAP):
\[
S_{\text{MAP}}\ldef (1-\alpha)I+\alpha P_{\vset}P_{\uset}
\]
with optimal $\alpha=\frac{2}{1+\sin(\theta_{F})^{2}}$, according
to~\cite{Bauschke_opt_rate_matr}.
\item Douglas-Rachford method (DR):
\[
S_{\text{DR}}\ldef \frac{1}{2}(I+R_{\vset}R_{\uset})
\]
 where $R_{\mathcal{C}}\ldef P_{C}^{2}=2P_{\mathcal{C}}-I.$
\item The optimal generalized alternating projections ($\text{GAP}^{*}$):
\[
S_{\text{GAP}^{*}}\ldef P_{\vset}^{\alpha^{*}}P_{\uset}^{\alpha^{*}},
\]
with $\alpha^{*}=\frac{2}{1+\sin\theta_{F}}.$
\item The adaptive generalized alternating projections (GAPA):
\[
S_{\text{GAPA}}\ldef P_{\vset}^{\alpha_{k}}P_{\uset}^{\alpha_{k}},
\]
implemented as in Algorithm~\ref{alg:adaptivegap}.
\item Generalized alternating projections with $a=1,\alpha_{1}=2,\alpha_{2}=\frac{2}{1+\sin(2\theta_{F})}$
(GAP$2\alpha$):
\[
S_{\text{GAP}2\alpha}=P_{\vset}^{\alpha_{2}}R_{\uset},
\]
as described in Section~\ref{sec:otherchoises}.
\item Generalized alternating projections with $\alpha=1,\alpha_{1}=\alpha_{2}=1.8$
(GAP1.8):
\[
S_{\text{GAP}1.8}\ldef P_{\vset}^{1.8}P_{\uset}^{1.8}.
\]
\end{itemize}
For each of the methods we monitor the shadow sequence
\[
z^{k}=P_{\uset}S^k x_{0}
\]
 and terminate when
\[
\norm{P_{\vset\cap\uset}z^{k}-z^{k}}<10^{-8}
\]
or when the number of iterations reach $200,000$.
\begin{rem}
The analysis in this paper concerns the convergence of the sequence towards a fixed-point.
We are actually more interested in the shadow sequence
(that we monitor in the examples),
since it can find a point in the intersection long before the
sequence converges to the fixed-point set.
This may be favorable for the Douglas-Rachford algorithm
because of its dominating complex eigenvalues,
compared to what its convergence rate suggests.
\end{rem}
The problems were generated and solved with \texttt{Julia}~\cite{julia}, and the results are shown in Figure~\ref{fig:iters-1-1}.
We see that the methods perform in line with the theoretical rates.
The method with optimal parameters performs considerably better and more reliably than for other choices.
We see that the adaptive method performs almost identically to the optimal parameters,
without prior knowledge of the Friedrichs angle.

We have verified numerically that the estimate in the
adaptive method converges to the Friedrichs angle.
For all problems that took more than 17 iterations to converge,
the estimate in the last iteration, was indeed conservative
($\hat\theta^k > \theta_F$).
Furthermore, the relative error $|\hat\theta^k - \theta_F|/\norm{\theta_F}$
was smaller than $5\%$ ($0.1\%$) at the last iteration, for all problems that ran more than
$100$ ($400$) iterations.
These results were obtained, even though no measures were taken to ensure $x^0\in\uset+\vset$.

\section{Conclusions}

We derived the optimal parameters for the generalized alternating projections
method for two linear subspaces. The optimal rate is considerably better than
previously analyzed parameters,
and we verify the results with an extensive set of numerical examples.
We also presented an adaptive method, that in practice is able to perform in line
with the optimal parameters, with no prior knowledge about the problem.

It remains as future work to study how the results apply to
more general feasibility problems.

\begin{appendix}
\section{Appendix}

\subsection{Proof of Theorem~\ref{thm:gapP1-1}}\label{app:gapP1-1}
We divide the proof into two cases: $p+q<n$ and $p+q\geq n$.
\begin{case}[$p+q<n$]\normalfont
For this case we can use the results in Proposition~\ref{prp:projections}.
For convenience of notation we introduce
\begin{align}
    f(\theta)&\ldef \frac{1}{2}\left(2-\alpha_{1}-\alpha_{2}+\alpha_{1}\alpha_{2}\cos^2\theta\right)\\
    g(\theta)&\ldef \sqrt{f(\theta)^{2}-(1-\alpha_{1})(1-\alpha_{2})}
\end{align}
so that (\ref{eq:eig12}) becomes $\lambda_{i}^{1,2} = f(\theta_i)\pm g(\theta_i)$.
For $\alpha_1=\alpha_2=\alpha^*=\frac{2}{1+\sin\theta_F}$ we get
$f(\theta_{F})=1-\alpha^*+{\alpha^*}^2c_F^2=\frac{1-\sin\theta_F}{1+\sin\theta_F}=\alpha^*-1$ and
$g(\theta_F)=0$.
Therefore, $\lambda_{F}^{1,2}=\alpha^{*}-1=\frac{1-\sin\theta_{F}}{1+\sin\theta_{F}}$.
We also see that $f(\pi/2)=1-\alpha^{*},\,g(\pi/2)=0$.
Since $f(\theta)$ is linear in $\cos^2\theta$ and $\left|f(\theta_{F})\right|=\left|f(\pi/2)\right|=\alpha^{*}-1$,
it follows that $\left|f(\theta_{i})\right|\leq\alpha^{*}-1$ for all
$\theta_{i}\in\left[\theta_{F},\pi/2\right]$. This means that the
corresponding $\lambda_{i}^{1,2}$ are complex with amplitudes
\begin{align*}
\left|\lambda_{i}^{1,2}\right| & =\sqrt{f(\theta_{i})^{2}+\left|f(\theta_{i})^{2}-(1-\alpha^{*})^{2}\right|}=\sqrt{(1-\alpha^{*})^{2}}\\
 & =\alpha^{*}-1\quad\forall i:\,\theta_{F}\leq\theta_{i}\leq\pi/2.
\end{align*}
Lastly, for $\theta_{i}<\theta_{F}$, we have $\theta_{i}=0$ and $\lambda_{i}^{1,2}=\left\{ 1,(1-\alpha^{*})^{2}\right\} $
as seen in Equation~\eqref{eq:T1matrix}. To conclude, we have

\begin{align*}
\left|\lambda_{F}^{1,2}\right| & =\left|\lambda^{3}\right|=\alpha^{*}-1\\
\left|\lambda^{4}\right| & =\left(\alpha^{*}-1\right)^{2}<\alpha^{*}-1\\
\left|\lambda_{i}^{1,2}\right| & =\alpha^{*}-1\quad\forall i:\,\theta_{F}\leq\theta_{i}\leq\pi/2\\
\left|\lambda_{i}^{1,2}\right| & =\left\{ 1,(1-\alpha^{*})^{2}\right\} \quad\forall i:\,\theta_{i}<\theta_{F},
\end{align*}
where $\lambda^3$ corresponds to eigenvalues in $1-\alpha_2$ and $1-\alpha_1$ in Proposition~\ref{prp:eigset}.
The eigenvalues $\lambda=1$ are semisimple since the matrix in~\eqref{eq:T1matrix}
is diagonal for $\theta_{i}=0$.
We therefore conclude, from
Fact~\ref{fct:limitexists} and \ref{fact:Convergent-at-rate},
that $\alpha_{1}=\alpha_{2}=\alpha^{*}$ results in that $T=S$ in Equation
(\ref{eq:GAP2}) is linearly convergent with any rate $\mu\in\left(\gamma^{*},1\right)$
where $\gamma^{*}=\alpha^{*}-1=\frac{1-\sin\theta_{F}}{1+\sin\theta_{F}}$
is a subdominant eigenvalue.
\end{case}

\begin{case}[$p+q\geq n$]\normalfont
    The following proof follows closely that in~\cite[p. 54]{Bauschke_lin_rate_Friedrich}.
    We can extend the space $\mathbb{R}^n$ with $k$ extra dimensions so that $p+q < n +k \rdef \bar n$.
    Let $\bar{\uset}\ldef\uset\times\{0_k\}$, $\bar{\vset}\ldef\vset\times\{0_k\}$, and therefore
    \begin{align*}
    P_{\bar\uset} =
        \begin{pmatrix}
    P_{\uset} & 0 \\
    0 & 0_k
    \end{pmatrix},\quad
    P_{\bar\vset} =
        \begin{pmatrix}
    P_{\vset} & 0 \\
    0 & 0_k
    \end{pmatrix}.
    \end{align*}
    It follows that
    \begin{align*}
    \bar S = \bar T \ldef P_{\bar\uset}^{\alpha^*}P_{\bar\vset}^{\alpha^*} = \begin{pmatrix}
    T & 0 \\
    0 & (1-\alpha^*)^2I_k
    \end{pmatrix},
    \end{align*}
    where $S=T=P_{\uset}^{\alpha^*}P_{\vset}^{\alpha^*}$.
    By using Case 1 on this matrix we conclude that
    $\gamma(\bar S)=\gamma^*$.
    From their definition, we see that the principal angles between $\uset$ and $\vset$ are the same as the ones between $\bar\uset$ and $\bar\vset$.
    Since $\bar T$ is block diagonal we have $\sigma(\bar T) = \sigma(T)\cup\{(1-\alpha^*)^2\}$. We therefore know that the eigenvalues $\lambda_F^{1,2}=\alpha^*-1$ in $\bar T$, corresponding to the Friedrichs angle, must also be in $T$. This means that $\gamma(S)\geq \alpha^*-1 = \gamma^*$.
    But $\norm{S^k-S^\infty}\leq\norm{\bar{S}^k-\bar{S}^\infty}$ so from Fact~\ref{fact:Convergent-at-rate} we know that $\gamma(S)\leq \gamma(\bar S)=\gamma^*$.
    We have therefore shown that $\gamma(S)=\gamma^*$, and the proof is complete.
\end{case}

\subsection{Lemmas}
\begin{lem}\label{lem:trdet}
The matrix $M:=(2-\alpha^{*})I+\frac{\alpha^{*}}{\alpha_{1}}(T_{1}^{i}-I)$
where $T_{1}^{i}$ is the matrix~\eqref{eq:T1matrix} corresponding to the angle $\theta_{F}$
has trace and determinant:
\begin{eqnarray*}
\text{tr}M & = & \frac{2}{(1+s)\alpha_{1}}\left(-\alpha_{1}-\alpha_{2}+\alpha_{2}\alpha_{1}c^{2}+2\alpha_{1}s\right)\\
\det M & = & \frac{4s(1-s)}{\alpha_{1}(1+s)^{2}}\left(-\alpha_{1}-\alpha_{2}+\alpha_{1}\alpha_{2}(1+s)\right),
\end{eqnarray*}
where $s:=\sin\theta_{F},\,c:=\cos\theta_{F}$.
\end{lem}

\begin{pf}
Let $s:=\sin\theta_{F},\,c:=\cos\theta_{F}$.
The matrix can be written
\begin{eqnarray*}
M & = & (2-\alpha^{*})I+\frac{\alpha^{*}}{\alpha_{1}}\left(\begin{pmatrix}1-\alpha_{1}s^{2} & \alpha_{1}cs\\
\alpha_{1}(1-\alpha_{2})cs & (1-\alpha_{2})(1-\alpha_{1}c^{2})
\end{pmatrix}-I\right)\\
 & = & \begin{pmatrix}2-\alpha^{*}-\alpha^{*}s^{2} & \alpha^{*}cs\\
\alpha^{*}(1-\alpha_{2})cs & 2-\alpha^{*}+\frac{\alpha^{*}}{\alpha_{1}}\left((1-\alpha_{2})(1-\alpha_{1}c^{2})-1\right)
\end{pmatrix}\\
 & = & \begin{pmatrix}2-\alpha^{*}(1+s^{2}) & \alpha^{*}cs\\
\alpha^{*}(1-\alpha_{2})cs & 2-\alpha^{*}+\frac{\alpha^{*}}{\alpha_{1}}\left(\alpha_{1}\alpha_{2}c^{2}-\alpha_{2}-\alpha_{1}c^{2}\right)
\end{pmatrix}.
\end{eqnarray*}
Using that $\alpha^{*}=\frac{2}{1+s}$, we can rewrite the diagonal elements
\[
2-\alpha^{*}(1+s^{2})=\alpha^{*}\left(1+s-(1+s^{2})\right)=\alpha^{*}s(1-s)
\]
and
\begin{align*}
2-&\alpha^{*}+\frac{\alpha^{*}}{\alpha_{1}}\left(\alpha_{1}\alpha_{2}c^{2}-\alpha_{2}-\alpha_{1}c^{2}\right)=\\
& = \alpha^{*}(1+s)-\alpha^{*}+\alpha^{*}\left(c^{2}(\alpha_{2}-1)-\frac{\alpha_{2}}{\alpha_{1}}\right)\\
 & = \alpha^{*}\left(s+c^{2}(\alpha_{2}-1)-\frac{\alpha_{2}}{\alpha_{1}}\right).
\end{align*}
 We can extract the factor $\alpha^{*}cs$ from the matrix and get
\[
M=\alpha^{*}cs\begin{pmatrix}\frac{1-s}{c} & 1\\
1-\alpha_{2} & \frac{s+c^{2}(\alpha_{2}-1)-\frac{\alpha_{2}}{\alpha_{1}}}{cs}
\end{pmatrix}.
\]
The trace is therefore given by
\begin{eqnarray*}
\text{tr}M & = & \alpha^{*}cs\left(\frac{1-s}{c}+\frac{s+c^{2}(\alpha_{2}-1)-\frac{\alpha_{2}}{\alpha_{1}}}{cs}\right)\\
 & = & \alpha^{*}\left(2s-s^{2}+c^{2}\alpha_{2}-c^{2}-\frac{\alpha_{2}}{\alpha_{1}}\right)\\
 & = & \frac{\alpha^{*}}{\alpha_{1}}\left(-\alpha_{1}-\alpha_{2}+\alpha_{2}\alpha_{1}c^{2}+2\alpha_{1}s\right)\\
 & = & \frac{2}{(1+s)\alpha_{1}}\left(-\alpha_{1}-\alpha_{2}+\alpha_{2}\alpha_{1}c^{2}+2\alpha_{1}s\right)
\end{eqnarray*}
and the determinant
\begin{eqnarray*}
\text{det}M & = & \left(\alpha^{*}cs\right)^{2}\left(\frac{\left(1-s\right)\left(s+c^{2}(\alpha_{2}-1)-\frac{\alpha_{2}}{\alpha_{1}}\right)}{c^{2}s}-\frac{\left(1-\alpha_{2}\right)c^{2}s}{c^{2}s}\right)\\
 & = & \alpha^{*2}s\left(\left(s+c^{2}(\alpha_{2}-1)-\frac{\alpha_{2}}{\alpha_{1}}-s^{2}-c^{2}s(\alpha_{2}-1)+s\frac{\alpha_{2}}{\alpha_{1}}\right)-\left(1-\alpha_{2}\right)c^{2}s\right)\\
 & = & \alpha^{*2}s\left(s+c^{2}(\alpha_{2}-1)-\frac{\alpha_{2}}{\alpha_{1}}-s^{2}+s\frac{\alpha_{2}}{\alpha_{1}}\right)\\
 & = & \alpha^{*2}s\left(s-1+\alpha_{2}c^{2}+\frac{\alpha_{2}}{\alpha_{1}}(s-1)\right)\\
 & = & \alpha^{*2}s(1-s)\left(-1+\alpha_{2}(1+s)-\frac{\alpha_{2}}{\alpha_{1}}\right)\\
 & = & \frac{\alpha^{*2}s(1-s)}{\alpha_{1}}\left(-\alpha_{1}-\alpha_{2}+\alpha_{1}\alpha_{2}(1+s)\right)\\
 & = & \frac{4s(1-s)}{\alpha_{1}(1+s)^{2}}\left(-\alpha_{1}-\alpha_{2}+\alpha_{1}\alpha_{2}(1+s)\right)
\end{eqnarray*}
\end{pf}

\begin{lem}\label{lem:realpart}
Under the assumptions $\alpha=\frac{\alpha^{*}}{\alpha_{1}}$, $\cmt{0\neq}\alpha_{1}\geq\alpha_{2}\cmt{\neq0}\note{>0}$
and $\theta_F<\pi/2$, then, the matrix $M$ in Lemma~\ref{lem:trdet} satisfies
\[
\left(\alpha_{1}\neq\alpha^{*}\text{ or }
\alpha_{2}\neq\alpha^{*}\right)
\Rightarrow\max\re\sigma(M)>0.
\]
\end{lem}
\begin{pf}
We show the equivalent claim
\[
\max\re\sigma(M)\leq0\Rightarrow\alpha_{1}=\alpha_{2}=\alpha^{*}.
\]

We have $\max\re\sigma(M)\leq0$ if and only if both eigenvalues of $M$
have negative or zero real part, which is equivalent to
\[
\lambda_{1}+\lambda_{2}\leq0\,\text{and}\,\lambda_{1}\lambda_{2}\geq0.
\]
which is equivalent to
\[
\text{tr}M\leq0\,\text{and}\,\text{det}M\geq0.
\]
Using Lemma~\ref{lem:trdet}, this can be written
\begin{eqnarray*}
\begin{cases}
\frac{2}{(1+s)\alpha_{1}}\left(-\alpha_{1}-\alpha_{2}+\alpha_{2}\alpha_{1}c^{2}+2\alpha_{1}s\right) & \leq 0\\
\frac{4s(1-s)}{\alpha_{1}(1+s)^{2}}\left(-\alpha_{1}-\alpha_{2}+\alpha_{1}\alpha_{2}(1+s)\right) & \geq 0
\end{cases},
\end{eqnarray*}
where $s\ldef\sin(\theta_F)$ and $c\ldef\cos(\theta_F)$.
Since \note{$\alpha_1>0$}, $s\in(0,1)$, this is equivalent to
\begin{subnumcases}{\label{eq:maxrem0}}
   \cmt{\sign{\alpha_1}(}\alpha_{1}+\alpha_{2}-\alpha_{2}\alpha_{1}c^{2}-2\alpha_{1}s\cmt{)} & $\geq 0$\label{eq:maxrem}\\
   \cmt{\sign{\alpha_1}(}-\alpha_{1}-\alpha_{2}+\alpha_{1}\alpha_{2}(1+s)\cmt{)} & $\geq 0$.\label{eq:maxrem2}
\end{subnumcases}
This implies that the sum is positive, i.e.
\begin{align*}
\cmt{\sign{\alpha_1}}\big(\alpha_{1}+\alpha_{2}-\alpha_{2}\alpha_{1}c^{2}&-2\alpha_{1}s\big)+
\cmt{\sign{\alpha_1}}\left(-\alpha_{1}-\alpha_{2}+\alpha_{1}\alpha_{2}(1+s)\right)\\
& = \cmt{\sign{\alpha_1}}(\alpha_{2}\alpha_{1}s^{2}-2\alpha_{1}s+\alpha_{1}\alpha_{2}s)\\
& = \cmt{|}\alpha_{1}\cmt{|}s\left(\alpha_{2}s-2+\alpha_{2}\right) \geq 0\\
\Leftrightarrow \quad & \alpha_{2}(1+s) \geq 2\\
\Leftrightarrow \quad & \alpha_{2} \geq \frac{2}{1+s}=\alpha^{*}.
\end{align*}
But then\cmt{, by assumption $\alpha_1\geq\alpha_2$ and thus $\sign{\alpha_1}=1$, and},~\eqref{eq:maxrem} implies
\[
\alpha_{1}+\alpha_{2}-\alpha_{2}\alpha_{1}c^{2}-2\alpha_{1}s\geq0
\]
\[
\Rightarrow(\text{since }\alpha_{2}\geq\alpha^{*})
\]
\[
\alpha_{1}+\alpha_{2}-\alpha^{*}\alpha_{1}c^{2}-2\alpha_{1}s\geq0
\]
which is equivalent to
\begin{eqnarray*}
\alpha_{1}+\alpha_{2}-\alpha^{*}\alpha_{1}c^{2}-2\alpha_{1}s & =\\
\alpha_{1}+\alpha_{2}-2\alpha_{1}(1-s)-2\alpha_{1}s & =\\
\alpha_{1}+\alpha_{2}-2\alpha_{1} & =\\
\alpha_{2}-\alpha_{1} & \geq & 0
\end{eqnarray*}
i.e.
\[
\alpha_{2}\geq\alpha_{1}.
\]
But by assumption $\alpha_{1}\geq\alpha_{2}$ so we know that \eqref{eq:maxrem0}
implies $\alpha_{1}=\alpha_{2}\geq\alpha^{*}$. Equation~\eqref{eq:maxrem}
yields
\begin{alignat*}{2}
&& \alpha_{1}+\alpha_{2}-\alpha_{2}\alpha_{1}c^{2}-2\alpha_{1}s &\geq0\\
&\Rightarrow \quad&
2\alpha_{1}-\alpha_{1}^{2}c^{2}-2\alpha_{1}s&\geq0\\
 &\Leftrightarrow &
2-\alpha_{1}c^{2}-2s&\geq0\\
&\Leftrightarrow & 2\frac{(1-s)}{c^{2}}&\geq\alpha_{1}\\
&\Leftrightarrow & \alpha^*=\frac{2}{(1+s)}&\geq\alpha_{1}.
\end{alignat*}
where the implication is from $\alpha_{1}=\alpha_{2}$.
We have therefore shown that $\alpha^{*}\geq\alpha_{1}=\alpha_{2}\geq\alpha^{*}$
i.e.
\begin{align*}
\max\re\sigma(M)\leq0\Rightarrow\alpha_{1}=\alpha_{2}=\alpha^{*}.
\end{align*}
\end{pf}

\subsection{Proof of Theorem~\ref{thm:gapiff}}\label{app:gapiff}
The first direction is proven by Theorem~\ref{thm:gapP1-1}.
To prove the other direction, we consider two cases.
\begin{case}[$p+q<n$]
\normalfont
Assume that $\alpha_{1}\geq\alpha_{2}$ and let
$\alpha=\hat{\alpha}\ldef \frac{\alpha^{*}}{\alpha_{1}}$.
The eigenvalues to the GAP operator $S$ in~\eqref{eq:GAP} are $1+\alpha\left(\lambda-1\right)$
where $\lambda$ are the eigenvalues to $T$.
Motivated by Proposition~\ref{prp:eigset}, we need to consider the following eigenvalues
\begin{align*}
\{ \lambda_{i}^{1,2}\} _{i\in1,\dots,p}\cap\left\{ 1-\alpha_{2}, 1-\alpha_{1},\,(1-\alpha_{2})(1-\alpha_{1})\right\}.
\end{align*}
For the eigenvalue $\lambda=1-\alpha_{1}$, we get
\begin{align}\label{eq:eigin1mina}
1+\hat{\alpha}(\lambda-1)=1+\frac{\alpha^{*}}{\alpha_{1}}(1-\alpha_{1}-1)=1-\alpha^{*}.
\end{align}

We now show that all choices of $\alpha,\alpha_{1},\alpha_{2}$ results
in an eigenvalue with real part larger than $\alpha^{*}-1$
unless $\alpha=1,\,\alpha_{1}=\alpha_{2}=\alpha^{*}$.

Consider the eigenvalues to $I+\hat{\alpha}(T_{1}^{i}-I)$ where $T_{1}^{i}$
is the matrix~\eqref{eq:T1matrix} corresponding to the angle $\theta_{F}.$
We have
\begin{equation}
\max\re\sigma\left(I+\hat{\alpha}(T_{1}^{i}-I)\right)>\alpha^{*}-1\label{eq:maxrealpart}
\end{equation}
if and only if
\begin{equation}
\max\re\sigma\left((2-\alpha^{*})I+\hat{\alpha}(T_{1}^{i}-I)\right)>0.\label{eq:eigpositive}
\end{equation}

By Lemma~\ref{lem:realpart} we know that~\eqref{eq:eigpositive} is true when $\alpha=\hat{\alpha}$, unless $\alpha_1=\alpha_2=\alpha^*$.
We therefore know that for $\alpha=\hat\alpha$, unless the optimal parameters are selected,
there will always be one eigenvalue of $S$ in
$1-\alpha^*$ and one, corresponding to $\theta_F$, with real part greater than $\alpha^*-1$.
We now consider the two cases $\alpha>\hat{\alpha}$ and $\alpha<\hat{\alpha}$. First note that $\alpha$ acts as a scaling of the eigenvalues relative to the point $1$, i.e. $(1-\alpha)+\alpha\lambda=1+\alpha(\lambda-1)$.
It is therefore clear that $\alpha>\hat\alpha$ will result in one eigenvalue with real part less than $1-\alpha^*=-\gamma^*$, and thus $\gamma(S)>\gamma^*$.

Similarly, any $\alpha<\hat\alpha$ will result in one eigenvalue ($\lambda_F^1$) with real part greater than $\alpha^*-1=\gamma^*$. If this eigenvalue is not in $1$, i.e. unless $1+\alpha(\lambda_F^1-1)=1$, we know that $\gamma(S)>\gamma^*$ also in this case.
Since $\alpha\neq0$ we have $1+\alpha(\lambda_F^1-1)=1$ if and only if $\lambda_F^1=1$.
But $\lambda_F^1=1$ only if $\det(T_F-I)=0$,
where $T_F$ is the block corresponding to $\theta_F$ in Equation~\eqref{eq:T1matrix}.
Since $\alpha_1,\alpha_2\neq0$ and $\theta_F>0$ we get
\begin{align*}
    \det(T_F-I)=
-\alpha_1s_F^2(\alpha_1c_F^2-\alpha_2+\alpha_1\alpha_2c_F^2)-\alpha_1^2(1-\alpha_2)c_F^2s_F^2=
\alpha_1\alpha_2s_F^2\neq0
\end{align*}
and thus $\lambda_F^1\neq1$.

We conclude that when $\alpha_1>\alpha_2$, then $\gamma(S)>\alpha^*-1$
for all parameters that are not $\alpha=1, \alpha_1=\alpha_2=\alpha^*$.

The proof is only dependent on the location of the eigenvalue corresponding to $\theta_F$ and the one in $1-\alpha_1$ from Proposition~\ref{prp:eigset}.
From symmetry of $\alpha_{1},\,\alpha_{2}$ in (\ref{eq:eig12}) we
see that the same argument holds if we instead assume $\alpha_{2}\geq\alpha_{1}$, let $\hat\alpha=\alpha^*/\alpha_2$, and pick the eigenvalue
$1-\alpha_{2}$, from the set in Proposition~\ref{prp:eigset}, in~\eqref{eq:eigin1mina} instead.
We have therefore shown that $\gamma(S)>\alpha^*-1$ unless $\alpha=1,\alpha_1=\alpha_2=\alpha^*$
so the result follows from Fact~\ref{fact:Convergent-at-rate}.
\end{case}

\begin{case}[$p+q\geq n$]
\normalfont
    As in the proof of Theorem~\ref{thm:gapP1-1} we can extend the space $\mathbb{R}^n$ with $k$ extra dimensions so that $p+q < n +k \rdef \bar n$.
    With $\bar{\uset}\ldef\uset\times\{0_k\}$, $\bar{\vset}\ldef\vset\times\{0_k\}$ 
    we get, as in the proof for Theorem~\ref{thm:gapP1-1},
    that $\sigma(\bar T)=\sigma(T)\cup\{(1-\alpha_1)(1-\alpha_2)\}$ and that the principal angles between $\uset$ and $\vset$ are the same as between $\bar\uset$ and $\bar\vset$.
    Using Proposition~\ref{prp:projections} on $\bar T$ we therefore see that the eigenvalues corresponding to the Friedrichs angle will exist in both $\bar T$ and $T$. The same is true for the eigenvalues in $1-\alpha_2$ and $1-\alpha_1$ for the respective cases $B1$ and $B3$ in Assumption~\ref{ass:dimUV}.
    The proof in Case 1 can therefore be used to show that $\gamma(S)>\gamma^*$ unless $\alpha=1, \alpha_1=\alpha_2=\alpha^*$.
\end{case}

\subsection{Proof of Theorem~\ref{thm:linconv}}\label{app:linconv}
Using~\cite[Thm. 2.12]{Bauschke_opt_rate_matr} we get for convergent $A$:
\begin{align*}
\norm{x^{k}-x^{*}} & =\norm{A^{k}x^{0}-A^{\infty}x^{0}}=\norm{(A^{k}-A^{\infty})x^{0}}\\
 & =\norm{(A-A^{\infty})^{k}x_{0}}\leq\norm{(A-A^{\infty})^{k}}\norm{x_{0}}.
\end{align*}
Using the spectral radius formula and $\rho(A-A^\infty)=\gamma(A)$~\cite[Thm. 2.12]{Bauschke_opt_rate_matr} we have, for any $\mu\in(\gamma(A),1)$
\begin{align*}
    \lim_{k\rightarrow\infty}\norm{(A-A^{\infty})^{k}}^\frac{1}{k} = \rho(A-A^\infty)=\gamma(A)<\mu,
\end{align*}
so there exists $N\in\mathbb{N}$ such that $\norm{(A-A^{\infty})^{k}}\leq\mu^k,\ \forall k\geq N$
and thus
\begin{align}
    \norm{x^{k}-x^{*}} \leq \mu^k\norm{x^0}\quad \forall k\geq N.
\end{align}

From~\cite[Corollary 2.7]{Bauschke_opt_rate_matr} we know that $S^{\infty}=P_{\fix S}$
since $S$ is nonexpansive, we therefore get $x^*=P_{\fix S}x^0$.

From Theorem~\ref{thm:gapP1-1} we know that $\gamma(S)=\frac{1-\sin\theta_F}{1+\sin\theta_F}$, and the proof is complete.

\end{appendix}

\end{document}